\documentstyle[twoside,12pt]{article}
\newcounter{beam}
\newenvironment{num}
{\begin{list}{(\roman{beam})}{\usecounter{beam}}
}{\end{list}}

\newenvironment{askn}
{\medskip\noindent {\bf Asknowledgements} \enspace
}{\medskip}

\newtheorem{utv}{Theorem}[section]
\newtheorem{lem}[utv]{Lemma}
\newtheorem{opr}[utv]{Definition}
\newtheorem{zam}[utv]{Remark}
\newtheorem{pri}[utv]{Example}
\newtheorem{sle}[utv]{Corollary}

\newtheorem{shag}{Step}

\def\proof{\par\noindent {\sl Proof:}\enspace}

\def\qed{{

\hfill{ Q.E.D.}}}

\def\kbi{{K+B }}
\def\kbn{{$K+B^{(n)}$ }}

\def\xni{{X^{(n)} }}
\def\xn{{$X^{(n)}$ }}
\def\xnsi{{X'^{(n)}}}
\def\xns{{$X'^{(n)}$ }}

\def\yni{{Y^{(n)} }}
\def\yn{{$Y^{(n)}$ }}
\def\zni{{Z^{(n)} }}

\def\bn{{$B^{(n)}$ }}
\def\bnsi{{B'^{(n)} }}
\def\bns{{$B'^{(n)}$ }}

\def\kbni{{K+B^{(n)}}{ }}
\def\xni{{X^{(n)}}{ }}
\def\bni{{B^{(n)}}{ }}

\title{Two Two-dimensional Terminations}
\author{Valery Alexeev}
\date{April 3, 1992}

\begin{document}
\maketitle


\def\alp{\alpha}		\def\Alp{\Alpha}
\def\bet{\beta}
\def\gam{\gamma}		\def\Gam{\Gamma}
\def\del{\delta}		\def\Del{\Delta}
\def\eps{\varepsilon}
\def\zet{\zeta}
\def\tet{\theta}		\def\Tet{\Theta}
\def\iot{\iota}
\def\kap{\kappa}
\def\lam{\lambda}		\def\Lam{\Lambda}
\def\sig{\sigma}		\def\Sig{\Sigma}
\def\vphi{\varphi}
\def\ome{\omega}		\def\Ome{\Omega}

\def\CalA{{\Cal A}}           \def\CalN{{\Cal N}}
\def\CalB{{\Cal B}}           \def\CalO{{\Cal O}}
\def\CalC{{\Cal C}}           \def\CalP{{\Cal P}}
\def\CalD{{\Cal D}}		 \def\CalQ{{\Cal Q}}
\def\CalE{{\Cal E}}           \def\CalR{{\Cal R}}
\def\CalF{{\Cal F}}           \def\CalS{{\Cal S}}
\def\CalG{{\Cal G}}           \def\CalT{{\Cal T}}
\def\CalH{{\Cal H}}           \def\CalU{{\Cal U}}
\def\CalI{{\Cal I}}           \def\CalV{{\Cal V}}
\def\CalJ{{\Cal J}}           \def\CalW{{\Cal W}}
\def\CalK{{\Cal K}}           \def\CalX{{\Cal X}}
\def\CalL{{\Cal L}}           \def\CalY{{\Cal Y}}
\def\CalM{{\Cal M}}           \def\CalZ{{\Cal Z}}

\def\bfa{{\bf a}}		\def\bfA{{\bf A}}
\def\bfb{{\bf b}}		\def\bfB{{\bf B}}
\def\bfc{{\bf c}}		\def\bfC{{\bf C}}
\def\bfd{{\bf d}}		\def\bfD{{\bf D}}
\def\bfe{{\bf e}}		\def\bfE{{\bf E}}
\def\bff{{\bf f}}		\def\bfF{{\bf F}}
\def\bfg{{\bf g}}		\def\bfG{{\bf G}}
\def\bfh{{\bf h}}		\def\bfH{{\bf H}}
\def\bfi{{\bf i}}		\def\bfI{{\bf I}}
\def\bfj{{\bf j}}		\def\bfJ{{\bf J}}
\def\bfk{{\bf k}}		\def\bfK{{\bf K}}
\def\bfl{{\bf l}}		\def\bfL{{\bf L}}
\def\bfm{{\bf m}}		\def\bfM{{\bf M}}
\def\bfn{{\bf n}}		\def\bfN{{\bf N}}
\def\bfo{{\bf o}}		\def\bfO{{\bf O}}
\def\bfp{{\bf p}}		\def\bfP{{\bf P}}
\def\bfq{{\bf q}}		\def\bfQ{{\bf Q}}
\def\bfr{{\bf r}}		\def\bfR{{\bf R}}
\def\bfs{{\bf s}}		\def\bfS{{\bf S}}
\def\bft{{\bf t}}		\def\bfT{{\bf T}}
\def\bfu{{\bf u}}		\def\bfU{{\bf U}}
\def\bfv{{\bf v}}		\def\bfV{{\bf V}}
\def\bfw{{\bf w}}		\def\bfW{{\bf W}}
\def\bfx{{\bf x}}		\def\bfX{{\bf X}}
\def\bfy{{\bf y}}		\def\bfY{{\bf Y}}
\def\bfz{{\bf z}}		\def\bfZ{{\bf Z}}

\def\Atil{{\widetilde A}}          \def\atil{{\tilde a}}
\def\Btil{{\widetilde B}}          \def\btil{{\tilde b}}
\def\Ctil{{\widetilde C}}          \def\ctil{{\tilde c}}
\def\Dtil{{\widetilde D}}          \def\dtil{{\tilde d}}
\def\Etil{{\widetilde E}}          \def\etil{{\tilde e}}
\def\Ftil{{\widetilde F}}          \def\ftil{{\tilde f}}
\def\Gtil{{\widetilde G}}          \def\gtil{{\tilde g}}
\def\Htil{{\widetilde H}}          \def\htil{{\tilde h}}
\def\Itil{{\widetilde I}}          \def\itil{{\tilde i}}
\def\Jtil{{\widetilde J}}      	   \def\jtil{{\tilde j}}
\def\Ktil{{\widetilde K}}          \def\ktil{{\tilde k}}
\def\Ltil{{\widetilde L}}          \def\ltil{{\tilde l}}
\def\Mtil{{\widetilde M}}          \def\mtil{{\tilde m}}
\def\Ntil{{\widetilde N}}          \def\ntil{{\tilde n}}
\def\Otil{{\widetilde O}}          \def\otil{{\tilde o}}
\def\Ptil{{\widetilde P}}          \def\ptil{{\tilde p}}
\def\Qtil{{\widetilde Q}}          \def\qtil{{\tilde q}}
\def\Rtil{{\widetilde R}}  	   \def\rtil{{\tilde r}}
\def\Stil{{\widetilde S}}	   \def\stil{{\tilde s}}
\def\Ttil{{\widetilde T}}          \def\ttil{{\tilde t}}
\def\Util{{\widetilde U}}          \def\util{{\tilde u}}
\def\Vtil{{\widetilde V}}          \def\vtil{{\tilde v}}
\def\Wtil{{\widetilde W}}          \def\wtil{{\tilde w}}
\def\Xtil{{\widetilde X}}          \def\xtil{{\tilde x}}
\def\Ytil{{\widetilde Y}}          \def\ytil{{\tilde y}}
\def\Ztil{{\widetilde Z}}          \def\ztil{{\tilde z}}

\def\abar{{\bar a}}	\def\Abar{{\bar A}}
\def\bbar{{\bar b}}	\def\Bbar{{\bar B}}
\def\cbar{{\bar c}}	\def\Cbar{{\bar C}}
\def\dbar{{\bar d}}	\def\Dbar{{\bar D}}
\def\ebar{{\bar e}}	\def\Ebar{{\bar E}}
\def\fbar{{\bar f}}	\def\Fbar{{\bar F}}
\def\gbar{{\bar g}}	\def\Gbar{{\bar G}}
\def\hBar{{\bar h}}	\def\Hbar{{\bar H}}
\def\ibar{{\bar i}}	\def\Ibar{{\bar I}}
\def\jbar{{\bar j}}	\def\Jbar{{\bar J}}
\def\kbar{{\bar k}}	\def\Kbar{{\bar K}}
\def\lbar{{\bar l}}	\def\Lbar{{\bar L}}
\def\mbar{{\bar m}}	\def\Mbar{{\bar M}}
\def\nbar{{\bar n}}	\def\Nbar{{\bar N}}
\def\obar{{\bar o}}	\def\Obar{{\bar O}}
\def\pbar{{\bar p}}	\def\Pbar{{\bar P}}
\def\qbar{{\bar q}}	\def\Qbar{{\bar Q}}
\def\rbar{{\bar r}}	\def\Rbar{{\bar R}}
\def\sbar{{\bar s}}	\def\Sbar{{\bar S}}
\def\tbar{{\bar t}}	\def\Tbar{{\bar T}}
\def\ubar{{\bar u}}	\def\Ubar{{\bar U}}
\def\vbar{{\bar v}}	\def\Vbar{{\bar V}}
\def\wbar{{\bar w}}	\def\Wbar{{\bar W}}
\def\xbar{{\bar x}}	\def\Xbar{{\bar X}}
\def\ybar{{\bar y}}	\def\Ybar{{\bar Y}}
\def\zbar{{\bar z}}  \def\Zbar{{\bar Z}}

\def\ahat{{\hat a}}		\def\Ahat{{\widehat A}}
\def\bhat{{\hat b}}		\def\Bhat{{\widehat B}}
\def\chat{{\hat c}}		\def\Chat{{\widehat C}}
\def\dhat{{\hat d}}		\def\Dhat{{\widehat D}}
\def\ehat{{\hat e}}		\def\Ehat{{\widehat E}}
\def\fhat{{\hat f}}		\def\Fhat{{\widehat F}}
\def\ghat{{\hat g}}		\def\Ghat{{\widehat G}}
\def\hhat{{\hat h}}		\def\Hhat{{\widehat H}}
\def\ihat{{\hat i}}		\def\Ihat{{\widehat I}}
\def\jhat{{\hat j}}		\def\Jhat{{\widehat J}}
\def\khat{{\hat k}}		\def\Khat{{\widehat K}}
\def\lhat{{\hat l}}		\def\Lhat{{\widehat L}}
\def\mhat{{\hat m}}		\def\Mhat{{\widehat M}}
\def\nhat{{\hat n}}		\def\Nhat{{\widehat N}}
\def\ohat{{\hat o}}		\def\Ohat{{\widehat O}}
\def\phat{{\hat p}}		\def\Phat{{\widehat P}}
\def\qhat{{\hapt q}}		\def\Qhat{{\widehat Q}}
\def\rhat{{\hat r}}		\def\Rhat{{\widehat R}}
\def\shat{{\hat s}}		\def\Shat{{\widehat S}}
\def\that{{\hat t}}		\def\That{{\widehat T}}
\def\uhat{{\hat u}}		\def\Uhat{{\widehat U}}
\def\vhat{{\hat v}}		\def\Vhat{{\widehat V}}
\def\what{{\hat w}}		\def\What{{\widehat W}}
\def\xhat{{\hat x}}		\def\Xhat{{\widehat X}}
\def\yhat{{\hat y}}		\def\Yhat{{\widehat Y}}
\def\zhat{{\hat z}}		\def\Zhat{{\widehat Z}}

\def\dbA{{\Bbb A}}       \def\dbN{{\Bbb N}}
\def\dbB{{\Bbb B}}       \def\dbO{{\Bbb O}}
\def\dbC{{\Bbb C}}       \def\dbP{{\Bbb P}}
\def\dbD{{\Bbb D}}       \def\dbQ{{\Bbb Q}}
\def\dbE{{\Bbb E}}       \def\dbR{{\Bbb R}}
\def\dbF{{\Bbb F}}       \def\dbS{{\Bbb S}}
\def\dbG{{\Bbb G}}       \def\dbT{{\Bbb T}}
\def\dbH{{\Bbb H}}       \def\dbU{{\Bbb U}}
\def\dbI{{\Bbb I}}       \def\dbV{{\Bbb V}}
\def\dbJ{{\Bbb J}}       \def\dbW{{\Bbb W}}
\def\dbK{{\Bbb K}}       \def\dbX{{\Bbb X}}
\def\dbL{{\Bbb L}}       \def\dbY{{\Bbb Y}}
\def\dbM{{\Bbb M}}       \def\dbZ{{\Bbb Z}}

\section{Introduction}

Varieties with log terminal and log canonical singularities are
considered in the Minimal Model Program, see \cite{kmm} for
introduction.  In \cite{sh:hyp} it was conjectured that many of the
interesting sets, associated with these varieties have something in
common: they satisfy the ascending chain condition, which means that
every increasing chain of elements terminates (in \cite{sh:hyp} it was
called the upper semi-discontinuaty). Philosophically, this is the
reason why two main hypotheses in the Minimal Model Program: existence
and termination of flips should be true and are possible to prove.

As for the latter, one of the main properties of flips is that log
discrepancies after doing one do not decrease and some of them
actually increase, \cite{sh:old}.  Therefore, if one could show that a
set of ``the minimal discrepancies'' satisfies the ascending chain
condition, that would help to prove the termination of flips.  The
Shokurov's proof of existence of 3-fold log flips \cite{sh:3f} is
another example of applying the same principle. In fact, to complete
the induction it uses some 1~-~dimensional statement, 2~-~dimensional
analog of which is proved in this paper. For further discussion, see
also
\cite{ag-kol}.

For one of the first examples where the phenomenon is actually proved
let us mention the following

\begin{utv}[\cite{al:fi},\cite{al:tg}]
Let us define the Gorenstein index of an $n$-dimensional Fano variety
$X$ with weak log terminal singularities as the maximal rational
number $r$ such that the anticanonical divisor $-K_X\equiv rH$ with an
ample Cartier divisor $H$. Then a set $$FS_n\cap
[n-2,=+\infty]=\{r(X)|X\; is\; a\; Fano\; variety\; and \;
r(X)>n-2\}$$ satisfies the ascending chain condition and has only the
following limit points: $n-2$ and $n-2+{1\over k}$, $k=1,2,3...$.
\end{utv}

In this paper we prove that the following two sets satisfy the
ascending chain condition:
\begin{num}
\item  (Theorems \ref{utv:local_pasc},\ref{utv:local_asc})
The set of minimal log discrepancies for $K_X+B$ where $X$ is a
surface with log canonical singularities.
\item (Theorem \ref{global_asc})
The set of groups $(b_1,...b_s)$ such that there is a surface $X$ with
log canonical and numerically trivial $K_X+\sum b_jB_j$.  The order on
such groups is defined in a natural way, see \ref{blessb'}.
\end{num}

\medskip

The proofs heavily use explicit formulae for log discrepancies from
\cite{al:lc}. We do not find it possible to prove them here again.
(This is quite easy anyway).

\begin{askn}
Author would like to thank V.V.Shokurov and J.Koll\'ar for asking the
questions that this paper gives the answers to and for useful
discussions.
\end{askn}

\section{Definitions and recalling}

All varieties in this paper are defined over the algebraically closed
field of characteristic zero.  $K_X$ or simply $K$ if the variety $X$
is clear from the context always denote the class of the canonical
divisor.

\subsection{Basics}

\begin{opr}
A {\bf \bfQ-divisor} on a variety $X$ is a formal combination $D=\sum
d_j D_j$ of Weil divisors with rational coefficients.
\end{opr}

\begin{opr}
One says that a \bfQ-divisor $D$ is {\bf \bfQ-Cartier} if some
multiple of it is a Weil divisor with integer coefficients that is a
Cartier divisor.
\end{opr}

\begin{opr}
\label{opr:discr}
Let $f:Y\to X$ be any birational morphism and $F_i$ be exceptional
divisors of this morphism. Consider a divisor of the form $K+B$, where
$B=\sum b_j B_j$ and $0<b_j\le1$. Coefficients $a_i$ in the following
formula $$K_Y+f^{-1}B+\sum F_i = f^*(K+B)+\sum a_i F_i$$ are called
{\bf log discrepancies} of $K+B$. \end{opr}

\begin{opr}
\label{opr:codiscr}
Let $f:Y\to X$ be any birational morphism and $F_i$ be exceptional
divisors of this morphism. Consider a divisor of the form $K+B$, where
$B=\sum b_j B_j$ and $0<b_j\le1$. Coefficients $b_i$ in the following
formula $$K_Y+f^{-1}B+\sum b_i F_i = f^*(K+B)$$ are called {\bf
codiscrepancies} of $K+B$. \end{opr}

\begin{zam}
Evidently there is a simple relation between log discrepancy and
codiscrepancy: $a_i=1-b_i$.
\end{zam}

\begin{opr}
A \bfQ-divisor of the form $K+B$ is said to be {\bf log canonical
(lc)} if
\begin{num}
\item $K+B$ is \bfQ-Cartier
\item there is a resolution of singularities $f:Y\to X$
 such that $supp(f^{-1}B)\bigcup F_i$ is a divisor with normal
intersections and all the log discrepancies $a_i\ge0$.  \end{num}
\end{opr}

\begin{opr}
A \bfQ-divisor of the form $K+B$ is said to be {\bf log terminal (lt)}
if
\begin{num}
\item $K+B$ is \bfQ-Cartier
\item there is a resolution of singularities
 $f:Y\to X$ such that $supp(f^{-1}B)\bigcup F_i$ is a divisor with
normal intersections and all the log discrepancies $a_i>0$.  \end{num}
\end{opr}

\subsection{Graphs}

With rare exceptions all the varieties in this paper will be
two-dimensional.  No doubt that the case of surfaces is much easier
than that of more-dimensional varieties. One of the reasons for this
is that surface has a natural quadratic form defined by intersection
of curves.  Many statements that we need can be formulated in terms of
weighted graphs and become therefore basicly combinatorical problems.

So let us start with a system of curves on a surface that are divided
into two classes: ``internal'', denoted by $F_i$ and ``external'',
denoted by $B_j$.

\begin{opr}
A weighted graph $\Gam$ is the following data:
\begin{num}
\item a ``ground graph'': each vertex $v$ of it corresponds to an ``internal''
curve
$F$ , two different vertices $v_1$ and $v_2$ are connected by wedge of
weight $F_1\cdot F_2$.
\item  weights: a vertex $v$ has   weight $w=-F^2$
\item  genera: a vertex $v$ has  genus $p_a(F)$ (arithmetical genus of
the curve)
\item  an ``external part'': additional vertices, corresponding to the
``external''
components $B_j$ , connected with vertices $v_i$ if $B_j$ and $F_i$
intersect.
\end{num}
\end{opr}

Vice versa, every weighted graph $\Gam$ corresponds to a system of
curves $\{F_i,B_j\}$.

\begin{opr}
Graph $\Gam$ is said to be {\bf elliptic, parabolic or hyperbolic} if
the corresponding quadratic form $F_i\cdot F_k$ is elliptic, parabolic
or hyperbolic, that is, has the signature $(0,n)$, $(0,n-1)$ or
$(1,n-1)$.
\end{opr}

The following is the basic case when we shall need such graphs: $X$ is
a surface with a divisor $K+B$ and $f:Y\to X$ is a resolution of
singularities of $X$.  The curves $F_i$ are exceptional divisors of
$f$ and the curves $B_j$ are strict transforms of the components of
$B$.  Note that since a matrix of intersection $(F_i\cdot F_k)$ is
negatively defined, the graph is elliptic and all the weights in this
case are positive integer numbers. Usually we will examine graphs that
correspond to the {\it minimal\/} resolution of singularities.

\begin{opr}
A graph $\Gam$ is said to be {\bf minimal} if it does not contain
internal vertices that have $p_a=0$ and weight 1.
\end{opr}

\begin{opr}
For any graph with a nondegenerate quadratic form $F_i\cdot F_k$ (for
example elliptic or hyperbolic) we define {\bf log discrepancies}
$a_i$ as the solutions of a system of linear equations $$\sum a_i
F_i\cdot F_k=(2p_a(F_k)-2-F^2)+(f^{-1}B+\sum F_i)F_k$$
\end{opr}

\begin{opr}
For any graph with a nondegenerate quadratic form $F_i\cdot F_k$ we
define {\bf codiscrepancies} $b_i$ by the formula $b_i=1-a_i$
\end{opr}

Let us explain the meaning of the two previous definitions. The
formulae above are equivalent to the following: $$(K+\sum b_jB_j+\sum
b_iF_i)F_k=0\quad for \quad any\quad k$$ So if the graph $\Gam$ is an
elliptic graph, corresponding to some birational morphism $f:Y\to X$
we get the previous definitions
\ref{opr:discr}, \ref{opr:codiscr}. Another situation when we shall
use discrepancies and codiscrepancies is the following: $X$ is a
surface with numerically trivial $K+B$, $f:Y\to X$ is some resolution.
Part of the vertices of $\Gam$ correspond to exceptional curves of $f$
and the other part -- to strict transforms of certain curves on $X$.

\begin{opr}
A graph $\Gam$ is said to be {\bf log canonical (lc) or log terminal
(lt)} with respect to $K+B$ if its log discrepancies $a_i\ge0$ or
$a_i>0$ respectively.
\end{opr}

 The main object into consideration in this paper is a surface $X$
with a divisor $K+B$ that is lc. So will be the corresponding graphs.
If we ignore the way $B$ meets the ground graph or assume that all the
coefficients of $B_j$ equal 1, then all such graphs are classified in
\cite{kaw:cb} (see also \cite{al:lc}).  They are divided into two
classes describing respectively rational and elliptic singularities.
In the case of rational singularities all the genera are equal to 0
(and by this reason will be omited), all the edges are simple (of
weight 1) the ground graphs are those of types $A_n$, $D_n$ and
$E_6,E_7,E_8$. If we fix some number $N$ and consider graphs with
weights $\le N$ then the only infinite series of such graphs are the
following

\begin{picture}(300,150)(-10,0)

\multiput(75,75)(150,0){2}{\begin{picture}(0,0)
\put(0,0){\oval(100,20)}
\put(-7.5,0){\oval(75,10)}
\put(40,0){\circle{10}}
\put(0,25){\vector(0,-1){15}}
\put(-15,-20){\vector(0,1){20}}
\put(-5,30){$m$}
\put(-20,-35){$q$}
\end{picture}}

\multiput(265,50)(0,50){2}{
\put(0,0){\circle{10}}
\put(10,-5){2}
}

\multiput(265,55)(0,25){2}{\put(0,0){\line(0,1){15}}}
\put(85,30){a)}
\put(235,30){b)}
\put(130,10){Figure 1}

\end{picture}

These are typical pictures that we shall use to describe graphs.  Long
ovals denote chains of vertices.  The numbers $q$ and $m$ denote the
absolute values of determinants of the submatrices of $F_i\cdot F_k$
that include only rows and columns corresponding to the vertices of
the chains.  It is very well known (\cite{ri}, comp.\cite{al:lc}) that
any chain is uniquely determined by a pair of coprime numbers $(q,m)$
with $1\le q<m$ and vice versa.  In the previous example $q$ and $m$
are any such numbers.

 Generally, graphs shall also have external parts that shall be
denoted by crossed vertices.

In the case of elliptic singularities one has ``circles'' of vertices
with $p_a=0$ and a single vertex with $p_a=1$.  $B$ is empty and all
the log discrepancies $a_i=0$, codiscrepancies $b_i=1$.

\begin{opr}
{\bf Du Val graph} is an elliptic graphs with all genera = 0, all
weights = 2 and empty external part $B$. It is well known that the
ground graph is then one of the graphs $A_n$, $D_n$, $E_6,E_7$ or
$E_8$.
\end{opr}

\begin{opr}
We say that a graph $\Gam'$ is a subgraph of $\Gam$ if all the
vertices of $\Gam'$ are at the same time vertices of $\Gam$, weights
of vertices and edges of $\Gam'$ and $p_a$ of vertices in $\Gam'$ do
not exceed the corresponding weights and $p_a$ in $\Gam$ and $F_i'\sum
b_j'B_j'\le F_i'\sum b_jB_j $ for the corresponding vertices.
\end{opr}

The following are easy linear algebra statements.

\begin{lem}
\label{a_less1}
Let $\Gam$ be a minimal elliptic graph.  Then all the log
discrepancies $a_i\le1$ (codiscrepancies $b_i\ge0$) and if $\Gam$ is
not a Du Val graph then $a_i<1$ ($b_i>0$).
\end{lem}
\proof Well known.
 \qed

\begin{lem}
\label{ag_less_ag'}
Let $\Gam'\subset\Gam$, $\Gam'\neq\Gam$ be two minimal elliptic graphs
and assume that the weights of the vertices are in both graphs the
same. Then for the log discrepancies one has $a_i\le a_i'$ (for
codiscrepancies $b_i\ge b_i'$) and if $\Gam$ is not a Du Val graph
then $a_i< a_i'$ ($b_i> b_i'$).  If the weights of $\Gam'$ and $\Gam$
are different then $a_i\le a_i'$ assuming that $\Gam$ is log
canonical.
\end{lem}
\proof Compare the corresponding systems of linear equations
(see \cite{al:fi}, \cite{al:lc}). \qed

\begin{lem}
\label{ag_greater_ag'}
Let $\Gam'\subset\Gam$, $\Gam'\neq\Gam$ be two graphs such that all
the log discrepancies of $\Gam'$ $a_i\le 1$ (codiscrepancies $b_i\ge
0$) and $v_0$ is a fixed vertex of $\Gam'$.  Assume that $\Gam$ is
hyperbolic and that $\Gam-v_0$ is elliptic.  Then for the log
discrepancy of $v_0$ one has $a_0\ge a_0' $ ($b_0\le b_0' $) assuming
that $\Gam$ is log canonical.
\end{lem}
\proof Compare the corresponding systems of linear equations. \qed

\begin{sle}
\label{max_ell}
Let $\Gam$ be a minimal elliptic graph and assume that all the log
discrepancies of $\Gam$ $a_i\ge c>0$. Then weights of the vertices are
bounded from above by $2/c$.
\end{sle}
\proof Consider a graph $\Gam'$ containing a single vertex of
weight $n$. Then $a'=2/n$.
\qed

\begin{sle}
\label{max_hyp}
Let $\Gam$ be a graph as in \ref{ag_greater_ag'} plus let $v_0$ have
weight 1.  Assume that the codiscrepancy of $v_0$ $b_0\ge c>0$. Then
$$\sum_{i\ne0} F_0F_i\le 2+\frac{2}{c}$$
\end{sle}
\proof Consider a graph $\Gam'$ containing a  vertex $v_0$ connected
with $n$ vertices of weight 2. Then $b_0'=2/(n+2)$.
\qed

\subsection{Sequences}

\begin{opr}
Let $X$ be a variety with a log canonical $K+B$.  {\bf A log
discrepancy of $\bf K+B$ $ld(K+B)$} is a minimal log discrepancy $a_i$
that appears in \ref{opr:discr} for some birational morphism $f:Y\to
X$.
\end{opr}

It is easy to see that $ld(K+B)$ is well defined and is a nonnegative
rational number.

\begin{opr}
Let $X$ be a surface with a log canonical $K+B$. {\bf A partial log
discrepancy of $\bf K+B$ $pld(K+B)$} is a minimal log discrepancy
$a_i$ that appears in
\ref{opr:discr} for the special  birational morphism $h:\Xtil\to X$, where
$\Xtil$ is the minimal resolution of singularities.
\end{opr}

\begin{opr}
Let $\xi=\{X^{(n)},K+B^{(n)}| n=1,2...\}$ be a sequence of surfaces.
Then we define $\bf ld(\xi)$ and $\bf pld(\xi)$ as the {\bf sequences}
of real numbers \{$ld(\kbni)$\} and \{$pld(\kbni)$\} respectively.
\end{opr}

\begin{opr}
\label{LDPLD}
Let $\xi=\{X^{(n)},K+B^{(n)}| n=1,2...\}$ be a sequence of surfaces.
Then we define $\bf LD(\xi)$ and $\bf PLD(\xi)$ as the {\bf subsets}
of real numbers \{$ld(\kbni)$\} and \{$pld(\kbni)$\} respectively.
\end{opr}

\begin{opr}
We define $ld,pld,LD,PLD$ for graphs in the same way as we have done it for
surfaces.
\end{opr}

\begin{opr}
\label{blessb'}
Let $B=(b_1,b_2...b_s)$ an $B'=(b_1',b_2'...b_t')$ be two groups of
numbers. One says that $B\le B'$ if

\begin{num}
\item $s\geq t$
\item for every $j=1...t$ $b_j\le b_j'$
\end{num}

If, in addition, one of the inequalities in (i) or in (ii) for some
index $j_0$ is strict, one says that $B<B'$.
\end{opr}

\begin{zam}
Because of the part (i) of \ref{blessb'} when considering a
nondecreasing sequence \bn we can always assume, passing to a
subsequence, that the lengths of \bn are in fact the same.
\end{zam}

\subsection{Log Del Pezzo surfaces}

\begin{opr}
A normal surface $X$ is said to be a Del Pezzo surface if $-K$ is an
ample \bfQ-divisor.
\end{opr}

The following is an simple lemma, see \cite{al-nik},\cite{nik} for the
proof which is especially easy if $K$ is lt or lc.
\begin{lem}
\label{Del_Pezzo}
Let $X$ be a log Del Pezzo surface an $h:\Xtil\to X$ be the minimal
resolution of singularities. Then
\begin{num}
\item the Kleiman-Mori cone  of effective curves $NE(\Xtil)$ is generated by
finitely many extremal rays
\item if $X\ne \bfP^2,\bfF_n$ (minimal rational surface) then all the extremal
rays are generated by exceptional curves of $f$ and (-1)-curves.
\end{num}
\end{lem}

\begin{lem}
\label{log_Del_Pezzo}
Let $X$ be a Del Pezzo surface and assume that $K$ is lc. Then $X$ is
one of the following:
\begin{num}
\item a rational surface with rational singularities
\item a generalized cone over a smooth elliptic curve
\end{num}
\end{lem}
\proof Let $h:\Xtil\to X$ be a minimal desingularization. $\Xtil$ is a
smooth surface and clearly $h^0(NK_{\Xtil}=0$ for any $N>0$, so
$\Xtil$ is ruled.

Assume that $X$ has a nonrational singularity. Then by the
classification of log canonical singularities $\Xtil$ contains an
elliptic curve or a circle of rational curves $F_0$ that is disjoint
from other curves, exceptional for $h$.  If $\Xtil$ is a locally
trivial $\bfP^1$-bundle then $F_0$ should be an exceptional section of
this bundle and should be smooth. In this case $X$ is a generalized
cone.  Otherwise $F_0$ should intersect a curve $E$ with $E^2<0$ that
lies in the fiber of a generically $\bfP^1$-bundle giving the
structure of a ruled surface and such that $E$ is not exceptional for
$h$.  By \ref{Del_Pezzo} $E$ is a $(-1)$-curve. The latter is
impossible since $-h^*K=-K_{\Xtil}-F_0-...$ and therefore $-h^*K\cdot
E\le0$.

Now let us assume that $X$ has only rational log canonical
singularities.  By the classification again one has
$-h^*K=-K_{\Xtil}-F_0-\sum b_iF_i$, $0\le b_i<1$ and $F_0$ is a
disjoint union of smooth rational curves.  Since $-h^*K$ is big, nef,
the Kawamata-Fiehweg vanishing gives
$$h^1(\Xtil,-F_0)=h^1(\Xtil,K_{\Xtil}+\sum b_iF_i+(-h^*K))=0$$ and
from the exact sequence $$0\to\cal O_{\Xtil}(-F_0)\to\cal
O_{\Xtil}\to\cal O_{F_0}\to 0$$ one gets $h^1(\cal O_{\Xtil})=0$.
Therefore $\Xtil$ and $X$ are rational surfaces.

\qed

\section{Local case: elliptic log canonical graphs}

{\bf In this section we consider only local situation. $\bfX$ is a
neighbourhood of a surface point $\bfP$ and all the components of
$\bfB$ pass through $\bfP$. }

\begin{utv}[Local boundness]
\label{utv:loc_bound}
Let $X,K+B$ be as above a neighbourhood of a surface point $P$ with lc
$K+B$ and all of $B_j$ pass through $P$. Then $\sum b_j\leq2$.
\end{utv}
\proof Proved in \cite{ag-kol} for the $n$-dimensional case with
a bound $n$.

\subsection{Minimal resolution}

\begin{utv}[Local partial ascending chain condition]
\label{utv:local_pasc}
Let $\xi=\{X^{(n)},K+B^{(n)}\}$ be a sequence of surfaces such that
\begin{num}
\item $K+B^{(n)}$ is lc
\item $B^{(n)}$ is a nondecreasing sequence
(for example, constant)
\end{num}
Then every increasing subsequence in $pld(\xi)$ terminates.

\noindent If, in the addition, one has
\begin{num}
\setcounter{beam}{2}
\item $B^{(n)}$ is an increasing sequence
\end{num}
then every nondecreasing subsequence in $pld(\xi)$ terminates.
\end{utv}

We prove \ref{utv:local_pasc} in several steps.

\begin{shag}
One can assume that \kbn and moreover $K$ are lt.
\end{shag}
\proof 	Indeed, if in (ii) \kbn is not lt, then
$pld(\kbni)=0$ but we are looking for increasing subsequences of
$pld(\xi)$.  In (iii) if \kbn is not lt, then there exists a partial
resolution $f:Y\to X$ with a single exceptional divisor $F$ such that
the corresponding log discrepancy $a=0$. Then
$$K_Y+F+f^{-1}B=f^*(K+B)$$ and the log adjunction formula for $F$ (see
\cite{sh:3f}) yields $$\sum {k-1+\sum l_j b_j\over k}=2$$ for some
positive integers $k,l_j$. It is easy to see that if \bn are
increasing then the sequence should terminate. \qed

\begin{shag}
By the previous step we can assume that there is a constant $\eps$ so
that for every $n$ $pld(\kbni)>\eps$.  Then we prove the following
\end{shag}

\begin{lem} All the lt elliptic graphs with $pld(\kbni)>\eps$
and $b_j>\eps$ can be described as follows:
\label{fin_many_graphs}
\begin{num}
\item finitely many graphs (that includes the way $B_j$ intersect $F_i$)
\item the graphs given on the next picture, where there are only finitely many
possibilities for the chains of vertices, denoted by ovals and for the
ways $B_j$ meet that vertices
\end{num}

\begin{picture}(320,270)(-10,-30)
\multiput(0,0)(0,100){2}{
\put(5,50){$q_1$}
\put(15,50){\vector(1,0){25}}
\put(65,50){\oval(80,20)}
\put(60,50){\oval(50,10)}
\put(95,50){\circle{5}}
\put(105,50){\line(1,0){12.5}}
\put(122.5,50){\line(1,0){15}}
\multiput(120,50)(20,0){2}{\circle{5}}
\multiput(117.5,55)(20,0){2}{2}
\multiput(142.5,50)(27.5,0){2}{\line(1,0){7.5}}
\multiput(155,50)(5,0){3}{\circle*{1}}
\multiput(35,85)(50,0){2}{\circle{14}}
\multiput(30,90)(50,0){2}{\line(1,-1){10}}
\multiput(30,80)(50,0){2}{\line(1,1){10}}
\multiput(50,80)(10,0){3}{\circle*{1}}
\put(35,78){\line(5,-6){15}}
\put(85,78){\line(-5,-6){15}}
\put(5,20){$m_1$}
\put(15,25){\vector(1,1){15}}
\put(125,20){\vector(-1,1){25}}
\put(125,15){$min$}
\put(110,85){\circle{14}}
\put(105,90){\line(1,-1){10}}
\put(105,80){\line(1,1){10}}
\put(95,52.5){\line(3,5){15}}
}
\multiput(180,30)(0,20){3}{\circle{5}}
\multiput(185,27.5)(0,20){3}{2}
\multiput(180,32.5)(0,20){2}{\line(0,1){15}}
\multiput(180,150)(20,0){2}{\circle{5}}
\multiput(177.5,155)(20,0){2}{2}
\put(182.5,150){\line(1,0){15}}
\put(202.5,150){\line(1,0){12.5}}
\put(255,150){\oval(80,20)}
\put(260,150){\oval(50,10)}
\put(225,150){\circle{5}}
\multiput(235,185)(50,0){2}{\circle{14}}
\multiput(230,190)(50,0){2}{\line(1,-1){10}}
\multiput(230,180)(50,0){2}{\line(1,1){10}}
\multiput(250,180)(10,0){3}{\circle*{1}}
\put(235,178){\line(5,-6){15}}
\put(285,178){\line(-5,-6){15}}
\put(210,185){\circle{14}}
\put(205,190){\line(1,-1){10}}
\put(205,180){\line(1,1){10}}
\put(225,152.5){\line(-3,5){15}}
\put(195,120){\vector(1,1){25}}
\put(175,115){$min$}
\put(305,150){\vector(-1,0){25}}
\put(310,150){$q_2$}
\put(310,120){$m_2$}
\put(305,125){\vector(-1,1){15}}
\put(140,-15){Figure 2}
\end{picture}

Moreover, the log discrepancies of any of suchs graphs satisfy the
following inequality $$pld(K+B)\ge {1-\sum l_jb_j\over m-q}$$
where $\bar b_j=\lim b_j$
 and tend
to this number as the chain of 2's gets longer and longer.  Here
$l_j=\sum( B_j\cdot F_i) r_i$, where $r_i$ is the determinant of the
short subchain of the ground graph, ``cut off'' by the vertex $v_i$.

\end{lem}

\proof By \ref{max_ell}  the weights of vertices in the graph $\Gam$ are
bounded
by $2/\eps$. Therefore, sacrificing finitely many graphs we can assume
that $\Gam$ is one of the graphs on Fig.1.

First, assume that we are in the case a) of Fig.1, i.e. $\Gam$ is a
chain.  Consider the sequence of log discrepancies of vertices in this
chain.  By \cite{al:lc} $$a_{i-1}-2a_i+a_{i+1}=(w_i-2)a_i+\sum
b_jB_jF_j\ge(w_i-2+\sum B_jF_i)\eps,$$ therefore the graph of this
function is concave up and, unless $w_i=2$ and all $B_jF_i=0$, it is
not a straight line but is ``very concave up''.  Now by \ref{a_less1}
the discrepancies $a_i\le1$. This implies that all the chains are
those on Fig.2 with only finitely many possibilities for the ovals and
with the chains of 2's of an arbitrary length $A$. Also, omiting
finitely many graphs, we can assume that the minimal log discrepancy
is achieved at one of the two vertices, where the arrows point out.

Now we use an explicit formula for the log discrepancies of those
vertices which follows easily from 3.1.8, 3.1.10 of \cite{al:lc}.

Define $\alp_1=1-\sum l_j^{(1)}b_j$ for the left part of the chain,
the meaning of $l_j$ being explained in the formulation of the
statement, $\alp_2=1-\sum l_j^{(2)}b_j$ be the corresponding
expression for the right par, and let $A$ be the length of the chain
of 2's. Then $$a_1={\alp_1(A(m_2-q_2)+m_2)+\alp_2q_1\over
A(m_2-q_2)(m_1-q_1)+m_2(m_1-q_1)+q_1(m_2-q_2)}$$ or
$$a_1={{\alp_1\over m_1-q_1}(A+{m_2\over m_2-q_2})+ {\alp_2\over
m_2-q_2}{q_1\over m_1-q_1}\over A+{m_2\over m_2-q_2}+{q_1\over
m_1-q_1}}$$ with the symmetric expression for $a_2$.

One can note that
\begin{enumerate}
\item if ${\alp_1\over m_1-q_1}\le {\alp_2\over m_2-q_2}$ then
${\alp_1\over m_1-q_1}\le a_1\le {\alp_2\over m_2-q_2}$
\item $\lim_{A\to\infty}a_1={\alp_1\over m_1-q_1}$
\end{enumerate}
and these two observations complete the proof in the case a) of Fig.1.
The case b) of Fig.1 is handled similarly. Let us mention only that in
the latter case there is only one possible vertex for the minimal log
discrepancy which is given by the formula $$a_1={\alp_1\over
m_1-q_1},$$ so this case can be treated formally as a subcase of a)
with $\alp_2=0$ and $m_2=q_2$.

\qed

\begin{shag}
The lemma~\ref{fin_many_graphs} implies \ref{utv:local_pasc}.
\end{shag}

\proof For any fixed graph $\Gam$ if the coefficients of the external part $B$
increase, then by \ref{ag_less_ag'} log discrepancy $pld(K+B)$
decreases.  Therefore, we can consider only case (ii) of
\ref{fin_many_graphs}.  Passing to a subsequence we can assume that
all the graphs are of the same type and the length of the sequence of
2's increases. But then $$pld(\kbni) \ge{1-\sum l_j\bar b_j\over
m-q}$$ and $$\lim pld(\kbni) = {1-\sum l_j\bar b_j\over m-q}$$ where
$\bar b_j=\lim b_j$, and we are done.

\qed

\begin{sle}\label{corol}
If $B=\emptyset$, then \ref{fin_many_graphs} says that the set of
minimal log discrepancies satisfies the ascending chain condition and
the only limit points are 0 and $1/k$, $k=2,3...$
\end{sle}

\begin{zam}
The statement \ref{corol} is due to V.V.Shokurov (unpublished).
\end{zam}

\subsection{General case}

Later we shall use the local ascending chain condition in the just
proved form, i.e.  for the minimal resolution of singularities. But
the minimal resolution of singularities of $X$ is not necessarily a
resolution of singularities for $K+B$, because $supp(f^{-1}B\cup F_i)$
can have nonnormal intersections.  Below we prove the statement,
corresponding to \ref{utv:local_pasc} but for $ld(\xi)$ instead of
$pld(\xi)$.  We first consider the case when \xn are nonsingular and
then combine our arguments to treat the general situation.

\begin{utv}
\label{nonsingular}
Let $\xi=\{X^{(n)},K+B^{(n)}\}$ be a sequence of {\bf nonsingular}
surfaces such that
\begin{num}
\item $K+B^{(n)}$ is lc
\item $B^{(n)}$ is a nondecreasing sequence
(for example, constant)
\end{num}
Then every increasing subsequence in $ld(\xi)$ terminates.

\noindent If, in the addition, one has
\begin{num}
\setcounter{beam}{2}
\item $B^{(n)}$ is an increasing sequence
\end{num}
then every nondecreasing subsequence in $ld(\xi)$ terminates.
\end{utv}
\proof
As above, we can assume that that \kbn are in fact lt.

Now let us find out what happens with a nonsingular surface $X$ with
$K+B$ after a single blow up $f:X\to Y$ at the point $P$, $F$ as
usually denotes the exceptional divisor of $f$. The answer is evident:
\begin{equation}
\label{blowup}
f^*(K+\sum b_jB_j)=K_Y+\sum b_jf^{-1}B_j+(-1+\sum mult_PB_jb_j)F
\end{equation}
and the condition $a>0$ translates to $-1+\sum mult_PB_jb_j<1$.  If
$-1+\sum mult_PB_jb_j\le0$, then for any further blowups all the log
discrepancies $a_i\ge a$, so they are irrelevant in finding the
minimal log discrepancy and $K+B$ is lt. However, if this is a
positive number, some negative log discrepancies can appear on the
following steps.

Now let $f:X\to Y$ be a composite of several blow ups. One gets

\begin{equation}
\label{blowups}
f^*(K+\sum b_jB_j)=K_Y+\sum b_jf^{-1}B_j+
\sum(-s_i+\sum t_{ik}b_k)F_i
\end{equation}
with some nonnegative integers $s_i$, $t_{ik}$ and $s_i\le\rho(Y/X)$.
The corresponding log discrepancies are given by $$a_i=1+s_i-\sum
t_{ik}b_k$$ and for fixed $s_i$ and nondecreasing/increasing $b_j$
they evidently form a nonincreasing/decreasing sequences. Note that
there are only finitely many such sequences with $a_i\ge0$.  Therefore
\ref{utv:local_pasc} follows from the following lemma.

\begin{lem}
With the assumptions as above, there is a constant $N(\xi)$ so that
for every surface \xn in $\xi$ there exists a birational morphism
$g:\yni\to\xni$ such that
\begin{num}
\item $\rho(\yni/\xni)<N(\xi)$
\item the minimal log discrepancy $ld(\kbni)$ is one of the log discrepancies
of $g$.
\end{num}
\end{lem}

\proof
Let us remind that we are in the local situation, so \xn is a
neighbourhood of a (nonsingular) point $P$.  Let $f:\zni\to\xni$ be a
single blow up at $P$. If in the formula~\ref{blowup} the number
$C=-1+\sum mult_PB_jb_j$ is positive and on \yn the strict transforms
of $B_j$ intersect at one point and have the same multiplicities as on
\xn, then by the formula~\ref{blowup} on the second blowup
codiscrepancy of the exceptional divisor equals $2C$, after the third
blowup $3C$ and so on (and it should be $\le1$).  Since $B^{(n)}$ is
nondecreasing, there exists a constant $\eps(\xi)$ so that for any
$-1+\sum m_jb_j>0$, one also has $-1+\sum m_jb_j>\eps(\xi)$. The
conclusion is that there exists a number $N_1$, depending on $\xi$, so
that after $N_1(\xi)$ blowups the configuration of $B_j$ simplifies in
some way: either the number of curves, passing through the points, or
the multiplicities at those points get smaller; or all the further
blowups are irrelevant in finding the minimal discrepancy.

Let $X^{(n)} {}'$ be $X^{(n)}$ with blown up points,
$\kbni'=f^*(\kbni)$. Note that the coefficients of $\bni'$ are still
nonnegative numbers. At the neighbourhood of any point of $X^{(n)}{
}'$\enspace $\bni'$ consists of several curves $B_j+\le 2$ nonsingular
curves $F_i$ with coefficients, given by the formula~\ref{blowups} and
hence, nondecreasing, and from the finite list of possible
combinations. Now we can find the next number $N_2(\xi)$ so that after
$N_2(\xi)$ blowups the configuration of ($B_j$ $+$ $\le 2$ nonsingular
curves) simplifies even further. By induction we get the desired
result.
\qed
\medskip

And finally we prove

\begin{utv}[Local  ascending chain condition]
\label{utv:local_asc}
Let $\xi=\{X^{(n)},K+B^{(n)}\}$ be a sequence of surfaces such that
\begin{num}
\item $K+B^{(n)}$ is lc
\item $B^{(n)}$ is a nondecreasing sequence
(for example, constant)
\end{num}
Then every increasing subsequence in $ld(\xi)$ terminates.

\noindent If, in the addition, one has
\begin{num}
\setcounter{beam}{2}
\item $B^{(n)}$ is an increasing sequence
\end{num}
then every nondecreasing subsequence in $ld(\xi)$ terminates.
\end{utv}

\proof By \ref{fin_many_graphs} all the singularities with $ld(\kbi)\ge\eps$
are divided into finite number of series $+$ finite number of graphs. The
latters are taken care by \ref{nonsingular}.

So all we have to do is to consider one of the graphs on Fig.2 with
the chain of 2's that is getting longer and longer.  And a simple
calculation shows that for all except finitely many graphs the minimal
log discrepancy is in fact one of log discrepancies of $h:\Xtil\to X$.
\qed

\section{Special hyperbolic log canonical graphs}

\noindent
{\bf Set-up\enspace } In this section $\Gam$ or $(X,K+B)$ always
denote the following:

\begin{opr} We say that a graph $\Gam$ is {\bf special hyperbolic} if
\begin{num}
\item $\Gam$ is hyperbolic and  connected
\item all the vertices have $p_a=0$, there is a special vertex
$v_0$ of weight 1, all other vertices $v_i$ have weights $\ge2$
\item $\Gam-v_0$ is elliptic
\item as usually, $\Gam$ may have an external part $B=\sum b_jB_j$
\end{num}
\end{opr}

Such graphs naturally appear when one considers a minimal resolution
of singularities of a Del Pezzo surface $X$ with $\rho(X)=1$ and $B_0$
being a (-1)-curve on the resolution.

\begin{utv}[Local-to-global ascending chain condition]
\label{hyper}
Let $\xi=\{\xni,\kbni\}$ be a sequence of special hyperbolic graphs
with a chosen vertex $B_0$ such that
\begin{num}
\item \kbn is lc
\item \{\bn\} is an increasing sequence, moreover, $\{b_0^{(n)}\}$ is an
increasing sequence
\item all the log discrepancies $a_i\ge 1-\bar b_0=1-\lim b_0$
\item $K+B^{(n)}$ is numerically trivial
\end{num}

Then $\xi$ terminates.
\end{utv}
\proof
{\par\noindent{\sl Case 1:}\enspace} $\bar b_0=\lim b_0=1$.

{}From \cite{al:lc} it follows that if $b_0$ is close enough to 1, then
all the singularities (that is, the connected components of
$\Gam-v_0$) and the ways the components of $B$ meet $F_i$ are
exhausted by the following list:

\begin{picture}(300,270)(-30,-30)
\multiput(0,0)(0,100){2}{
\put(20,45){$B_0$}
\put(50,50){\circle{14}}
\put(45,55){\line(1,-1){10}}
\put(45,45){\line(1,1){10}}
\put(130,50){\oval(110,20)}
\put(57,50){\line(1,0){18}}}
\multiput(170,20)(0,30){3}{\circle{10}}
\multiput(170,25)(0,30){2}{\line(0,1){20}}
\multiput(180,17)(0,60){2}{2}
\put(85,150){\circle{5}}
\put(135,150){\oval(80,10)}
\put(115,180){\vector(0,-1){25}}
\put(105,180){$q$}
\put(145,120){\vector(0,1){25}}
\put(150,115){$m$}
\multiput(135,185)(40,0){2}{\circle{14}}
\multiput(130,190)(40,0){2}{\line(1,-1){10}}
\multiput(130,180)(40,0){2}{\line(1,1){10}}
\multiput(150,180)(5,0){3}{\circle*{1}}
\put(135,178){\line(5,-6){15}}
\put(175,178){\line(-5,-6){15}}
\put(110,-15){Figure 3}
\end{picture}

The next step is a formula for the coefficient $b_0$ that follows from
the explicit calculations of \cite{al:lc}:

\begin{equation}
\label{horrible}
b_0={\sum_{s=1}^{N}{{q_s+\alp_s}\over m_s}-(N-1)\over
\sum_{s=1}^{N} {q_s\over
m_s}-1}=1-{\frac{(N-2)-\sum_{s=1}^{N}{\frac{\alp_s}{m_s}}}
{\sum_{s=1}^{N}{\frac{q_s}{m_s}}-1}}
\end{equation}
with denominator $>0$, where $N$ is a number of connected components
of $\Gam-v_0$.  Here $\alp_s=1-\sum l_j^sb_j$ as in
\ref{fin_many_graphs}.  We consider the second case of the figure~3
formally as a subcase of the first one with $q=m$ and $\alp=0$.

Note that by \ref{max_hyp} a number of graphs that $v_0^{(n)}$ is
connected with in the sequence is bounded.

For any fixed $N$ the conditions $\lim b_0=1$ and $b_0<1$ imply

$$ \sum_{s=1}^N{\frac{\alp_s}{m_s}}<N-2\quad and\quad
\lim \sum_{s=1}^N{\frac{\alp_s}{m_s}}=N-2.$$

We can assume that some of $m_s$ are fixed and others tend to
infinity. For the latters ${\frac{\alp_s}{m_s}}\to 0$ and
${\frac{\alp_s}{m_s}}>0$.  This is so by \ref{fin_many_graphs} (here
it is important again that there is a constant $\eps(\xi)$ so that
$\sum m_jb_j-1>0$ implies $\sum m_jb_j-1>\eps(\xi)$) and by
\ref{utv:loc_bound}. So we can assume that
$\sum_{s=1}^M{\frac{\alp_s}{m_s}}<N-2$ and $$\lim
\sum_{s=1}^M{\frac{\alp_s}{m_s}} =\lim\sum_{s=1}^M{\frac{1-\sum
l_j^sb_j}{m_s}}=N-2$$ with $m_1...m_M$ being fixed. But this
definitely gives a contradiction. Note that $\sum l_j^sb_j\le 2$ by
\ref{utv:loc_bound}, so there are only finitely many possibilities for
$l_j^s$.

Finally, for $N\ge5$ $$b_0\le1-{\frac{(N-2)-\sum{1\over{m_s}}}{N-1}}
\le 1-{\frac{{N\over 2}-2}{N-1}}\le {7\over 8}$$
\noindent and we are done.

\medskip
{\par\noindent{\sl Case 2:}\enspace} $\bar b_0=\lim b_0<1$.

Since all the log discrepancies $a_i\ge\eps=1-{\bar b}_0$, the only
infinite series of connected components of $\Gam-v_0$ are given by
\ref{fin_many_graphs}. Moreover, for the minimal log discrepancies
there $$\lim\min a_i\le{\frac{1-\sum(B_0F_i)r_ib_0}{m-k}}$$

\noindent and this should be not less than $1-\bar b_0$.
As a conclusion, all the infinite series are given by

\begin{picture}(300,250)(-50,-40)
\multiput(0,0)(0,100){2}{
\put(-5,45){$B_0$}
\put(20,50){\circle{14}}
\put(15,55){\line(1,-1){10}}
\put(15,45){\line(1,1){10}}
\put(27,50){\line(1,0){15.5}}
\multiput(45,50)(20,0){2}{\circle{5}}
\multiput(42.5,55)(20,0){2}{2}
\put(47.5,50){\line(1,0){15}}
\multiput(67.5,50)(27.5,0){2}{\line(1,0){7.5}}
\multiput(80,50)(5,0){3}{\circle*{1}}}

\multiput(105,150)(20,0){2}{\circle{5}}
\multiput(102.5,155)(20,0){2}{2}
\put(107.5,150){\line(1,0){15}}

\put(127.5,150){\line(1,0){17.5}}
\put(190,150){\oval(90,20)}
\put(155,150){\circle{5}}
\put(195,150){\oval(60,10)}
\multiput(155,185)(70,0){2}{\circle{14}}

\multiput(150,190)(70,0){2}{\line(1,-1){10}}
\multiput(150,180)(70,0){2}{\line(1,1){10}}
\multiput(180,180)(10,0){3}{\circle*{1}}
\put(155,178){\line(5,-6){15}}
\put(225,178){\line(-5,-6){15}}

\multiput(105,30)(0,20){3}{\circle{5}}
\multiput(110,28.5)(0,20){3}{2}
\multiput(105,32.5)(0,20){2}{\line(0,1){15}}
\put(125,-15){Figure 4}
\end{picture}

Now we would like to use a variant of the formula~\ref{horrible}.
However, $B_0$ can intersect finitely many types of graphs
arbitrarily. Still, for any fixed combination, if $b_j$ increase,
$b_0$ decreases. The situation is exactly the opposite to the one of
elliptic graphs since the signature of the quadratic form is now
$(1,n-1)$ instead of $(0,n)$ and the graph $\Gam-v_0$ is still
elliptic (cf. \ref{max_hyp}).

All the said above implies that for $b_0$ there are only finitely many
possible expressions of the form $$b_0=1-{C_1+\sum
C_2^jb_j-\sum_{s=1}^{N}{\alp_s\over m_s}\over
C_3+\sum_{s=1}^{N}{q_s\over m_s}}$$ with fixed
$C_1,C_2^j,C_3,m_s-q_s$, $m_s\to+\infty$, $C_2^j\ge0$ and the
denominator $>0$.  Simplifying, $$1-b_0={C_1+\sum
C_2^jb_j-\sum_{s=1}^{N}{\alp_s\over m_s}\over
C_3'-\sum_{s=1}^{N}{m_s-q_s\over m_s}}$$

Now $\lim b_0=\bar b_0$ implies $(C_1+\sum C_2^j\bar b_j)/C_3'=1-\bar
b_0$.  And, finally, the inequalities ${\alp_s\over m_s-q_s}\ge 1-\bar
b_0$ and $C_2^j\ge0$ imply that$1-b_0\le1-\bar b_0$, i.e. $b_0\ge \bar
b_0$ that gives a contradiction.

\qed

\begin{zam}
As the proof shows,
\ref{hyper} is not true without the assumption (iii).
\end{zam}

\section{Global case}

\begin{utv}[Global boundness]
\label{utv:glob_bound}
Let $X$ be a surface with a lc divisor $K+B$ and assume that $f:X\to
Y$ is a contraction of an extremal ray such that $K+B$ if
$f$-nonpositive.  Let $B^+=\sum b_j^+B_j^+$ contain all the components
in $B$ that are $f$-positive.  Then
\begin{num}
\item if $\dim Y=2$, $\sum b_j^+\le2$
\item if $\dim Y=1$, $\sum b_j^+\le2$
\item if $\dim Y=0$, $\sum b_j^+\le3$
\end{num}
\end{utv}
\proof
(i) follows from \ref{utv:loc_bound}, because $K_Y+f(B)$ is also lc.
(ii) is clear: if $B^+$ is not empty, then $-K$ should be negative on
the general fiber, so a general fiber $F$ is isomorphic to $\bfP^1$
and $\sum b_j^+\le B^+F\le -KF=2$.

In the case (iii) if $X$ is nonsingular, then $X\simeq {\bfP^2}$ and
the statement is evident. If $X$ is singular, consider a partial
resolution $g:Z\to X$, dominated by the minimal resolution and such
that $\rho(Z)=\rho(X)+1=2$. Then by \ref{Del_Pezzo} there is a second
extremal ray and $g^*(K+B)=K_Z+B_Z$ is nonpositive with respect to
this extremal ray. Since every curve on $Z$ is positive with respect
to at least one of the extremal rays, (iii) with the bound 4 follows
immediately. To get the bound 3 is an easy exercise.

\qed

\begin{zam}
\ref{utv:glob_bound}(iii)
 is also proved in \cite{ag-kol} for arbitrary dimension with a bound
$n+1$.
\end{zam}

\begin{utv}[Global  ascending chain condition]
\label{global_asc}
Let $\xi=\{X^{(n)},K+B^{(n)}\}$ be a sequence of surfaces such that
\begin{num}
\item $K+B^{(n)}$ is lc
\item $B^{(n)}$ is an increasing sequence
\item $K+B^{(n)}$ is numerically trivial
\end{num}
Then $\xi$ terminates.
\end{utv}

Proof of \ref{global_asc} will be given in several steps.

\setcounter{shag}{0}
\begin{shag}
One can assume that all the surfaces \xn are Del Pezzo surfaces with
$\rho(\xni)=1$.
\end{shag}
\proof
We can assume that the lengths of the groups \bn in the sequence $\xi$
are constant and that $b_1$ always increases.  Now consider a divisor
$K+B-\eps B_1$ on $X^{(n)}$.
Note here that $B_1$ is $\bfQ$-factorial by the classification of log
canonical singgularities.
 It is lc and is not numerically effective
and if $B_1^2\le0$ then $(K+B-\eps B_1)B_1\ge0$.  Therefore either
$\rho(\xni)=1$ and then \xn is a Del Pezzo surface with lc $K+B$ or
there is an extremal ray that does not contract $B_1$.  If the
contraction is birational, we make it and repeat the same procedure
again. If it is a fibration, the claim follows from the corresponding
1-dimensional statement.

\qed

\begin{zam} The argument works in the 3-dimensional case as well.
\end{zam}

\begin{shag}
One can assume that there are only finitely many different types of
graphs of singularities that the increasing components of \bn are
passing through.
\end{shag}
\proof
As usually, we can assume that the groups \bn have the same length.
Now consider the set $PLD(\xi)$. By \ref{utv:local_pasc} this set
satisfies the ascending chain condition and has at least one limit
point. Let $l$ be the minimal limit point of $PLD(\xi)$. Fix the
number $C$ so that all $b_j\ge C$. If the surfaces in $\xi$ contain
singularities that correspond to infinitely many elliptic graphs, then
by \ref{fin_many_graphs} $l\le1-C$. Passing to a subsequence we can
assume that a sequence of minimal log discrepancies, which we shall
denote $\{a^{(n)}_s\}$ is a decreasing sequence and $\lim a^{(n)_s}=l$
(the sequence of codiscrepancies is increasing and $\lim
b_0^{(n)_s}=1-l\ge C$.

Now consider a partial resolution $f:\yni\to\xni$ which is dominated
by the minimal desingularization and which blows up exactly the curve
$B^{(n)_s}$. Then $$f^*(K+\bni)=K_Y+f^{-1}\bni+b^{(n)_s}B^{(n)_s}$$

The surface \yn has Picard number 2 and by \ref{Del_Pezzo} there is a
second extremal ray, corresponding to a $(-1)$-curve on $\Ytil=\Xtil$.
Let $g:\yni\to\xnsi$ be the contraction of this second extremal ray.
If $g$ is a fibration then restricting of
$K_Y+f^{-1}\bni+b^{(n)_s}B^{(n)_s}$ on the general fibre of $g$ readily
gives a contradiction. Hence, we shall assume that $g$ is a birational
morphism.  A divisor $K+\bnsi=g_*f^*(k+\bni)$ is lc and numerically
trivial, \bns has either the same number of components as \bn or one
more, and, after passing to a subsequence, $B'^{(n)}$ is an increasing
sequence.

A morphism $g$ can contract one of the components of \bn and we can
assume that it is always, say, $B_0$.  However, by \ref{hyper} and
\ref {utv:local_pasc} the sequence $\{b_0^{(n)}\}$ cannot be an
increasing sequence with $\lim b_0^{(n)}\ge1-l$. Therefore, changing the
sequence $\xi=\{X^{(n)}\}$ by a new sequence $\xi'=\{X'^{(n)}\}$, we
are gaining a new component with increasing coefficient that has the
limit $1-l\ge C$.  Note that for a new minimal limit point $l'$ of
$PLD(\xi')$ one has $l'\ge l$.  This is so because a minimal
desingularization of \xns is dominated by the minimal
desingularization of \xn and $K+B^{(n)}$, $K+B'^{(n)}$ both are
numerically trivial, so $PLD(\xi')$ is a subset in $PLD(\xi)$.

Repeating the procedure, we get one more component and so on. After
$k$ steps the sum of the coefficients in \bn will be greater than
$kC$. This eventually will get into the contradiction with
\ref{utv:glob_bound}.

\qed

\begin{shag}
One can assume that all the surfaces \xn are isomorphic to each other.
\end{shag}
\proof
By \ref{log_Del_Pezzo} a surface $\Xtil^{(n)}$ is either a locally
trivial $\bfP^1$-bundle with a section which is a smooth elliptic
curve or a rational surface with rational singularities. In the former
case the statement follows from the 1-dimensional analog by
restricting $B$ to the fiber of the fibration.  Now assume we are in
the latter case.  By the previous step, there exists a constant
$N(\xi)$ so that for the increasing component $B_1$ of $B$\enspace
$NB_1$ is Cartier. Hence for any curve $D$ on
\xn
$$-KD=\sum b_jB_jD\ge b_1/N\ge C/N$$ Now theorem~2.~3 of \cite{al:fi}
states that for all such surfaces $\rho(\Xtil)$ is bounded.  Therefore
one can get $X^{(n)}$ by blowing up finitely many points from the
minimal rational surface $\bfF_k$. Threfore there are only finitely
many possibilities for the graph of exceptional curves on $X^{(n)}$
except for the fact that one weight $k$ can be arbitrary.  Now if
$B_1$ does not lie in the fiber for infinitely many $n$ we prove the
statement restricting a numerically trivial divisor $K+\sum
b_jB_j+\sum b_i F_i$ to the fiber and using \ref{utv:local_pasc}.
Otherwise (recall that $\rho(X)=1$) $B_1$ on $X$ should pass through
the singularity which graph contains the exceptional curve of
$\bfF_k$. By the previous step $k$ is bounded. Hence we can assume
that the surfaces $X^{(n)}$ belong to a bounded family and it is
enough to consider only finitely many of them.

\qed

\begin{zam}
Theorem 2.3 in \cite{al:fi} is stated for log terminal singularities.
But in fact the proof is exactly the same for rational log canonical
singularities.
\end{zam}

\begin{shag}
\ref{global_asc} follows.
\end{shag}
\proof
Indeed, there are only finitely many possibilities for effective Weil
divisors $B_j$.

\qed
\medskip

The following example shows that \ref{global_asc} is not true without
the assumption~(i).

\begin{pri}
Consider a sequence of surfaces \xn so that $\Xtil^{(n)}=\bfF_n$ and
$\bni=(1-1/n)F_1+3/4(F_2+F_3+F_4)$ where $F_1$ is an image of the
infinite section of $\bfF_n$, $F_{2,3,4}$ are fibres. Note that \kbn
is numerically trivial but is not lc.
\end{pri}

\end{document}